\documentclass[10pt]{amsart}

\usepackage{amsmath}
\usepackage{amssymb}
\usepackage{amscd}
\usepackage{epsfig}
\usepackage{color}
\usepackage{subfigure}
\usepackage{graphicx}

 \newcommand{\bee}{\begin{equation}}
 \newcommand{\eee}{\end{equation}}
 
 \newcommand{\Lb}{\mbox {\boldmath ${\Lambda}$}}
 \newcommand{\Gb}{\mbox {\boldmath ${\Gamma}$}}
 
 \newcommand{\Lbs}{\mbox{\scriptsize\boldmath ${\Lambda}$}}

\def\b0{{\bf 0}}

\definecolor{darkgreen}{rgb}{0.3,0.6,0.1}

\newcommand{\be}{\begin{eqnarray}}
\newcommand{\ee}{\end{eqnarray}}
\newcommand{\supp}{\mbox{\rm supp}}

\newcommand{\freq}{\mbox{\rm freq}}
\newcommand{\Vol}{\mbox{\rm Vol}}

\newcommand{\R}{{\mathbb R}}

\newcommand{\Z}{{\mathbb Z}}
\newcommand{\C}{{\mathbb C}}

\newcommand{\Ak}{{\mathcal A}}

\newcommand{\Dk}{{\mathcal D}}

\newcommand{\Sk}{{\mathcal S}}
\newcommand{\Tk}{{\mathcal T}}

\newcommand{\Lam}{{\Lambda}}

 \newtheorem{theorem}{Theorem}[section]

 \newtheorem{defi}[theorem]{\it{Definition}}
\numberwithin{equation}{section}

\begin{document}

\title{The computation of overlap coincidence in Taylor-Socolar substitution tiling }


\date{August 24, 2012}

\thanks{*Corresponding author. \\
This research was supported by Basic Research Program through the
National Research Foundation of Korea(NRF) funded by the Ministry
of Education, Science and Technology(2010-0011150) and the
Japanese Society for the Promotion of Science (JSPS), Grant in aid
21540010. The second author is grateful for the support of KIAS for this research.   }

\maketitle

 \centerline{
{ Shigeki Akiyama $^{\,\rm a}$ and Jeong-Yup Lee $^{\,\rm b*}$}}

\hspace*{4em}

\smallskip

\begin{abstract}

Recently Taylor and Socolar introduced an aperiodic mono-tile. 
The associated tiling can be viewed as a substitution tiling. 
We use the substitution rule for this tiling
and apply the algorithm of \cite{AL} to check overlap coincidence.
It turns out that the tiling has overlap coincidence. So the
tiling dynamics has pure point spectrum and we can conclude that this
tiling has a quasicrystalline structure.

\hspace*{4em}


\end{abstract}

\section{Introduction}

{\it Aperiodic tiles} are the set of prototiles 
which tile the space with their isomorphic images by Euclidean motions
(composition of translations, rotations and reflections)
but only in non-periodic way. There have been many examples of aperiodic tiles and study on them \cite{AGS, Berg, Robinson, Gardner, GS, Morita, MT, Radin, Senechal}.
Two of well-known examples of aperiodic tiles with simple prototiles up to Euclidean motions are
Penrose tiles and Ammann tiles which are uncovered in the mid 70's. These sets consist of 
two prototiles and it has been the smallest number of prototiles which form aperiodic 
tiles until recently. 
Since then, people have been interested in finding a single prototile for an aperiodic tile. This problem is coined as `Mono-tile' problem or `Einstein' problem (one stone 
in German).
It had taken quite some time before Taylor and Socolar announced in 2010 the existence of an aperiodic mono-tile.
Their tile is a hexagonal tile with colored decorations and matching rules which can be embedded onto a single tile using shape only. Penrose had found earlier a mono-tile using matching rules which is the reformulation of $(1 + \epsilon + \epsilon^2)$ aperiodic tiles given in \cite{Penrose}, but in this case it cannot be replaced by a single tile using only the shape. One needs two other tiles to replace the matching rules in the Penrose mono-tile. Both of Taylor-Socolar mono-tile and Penrose functional mono-tile are based on hexagonal shape. But they make different tilings. 
It would be interesting to understand the tiling space of Penrose mono-tile tilings and compare with the tiling space of Taylor-Socolar tilings. 

Around the mid 80's, Shechtman \cite{She} discovered a quasicrystal
with forbidden rotational symmetry of crystal diffraction pattern. After the confirmation of the existence of quasicrystals, crystal is redefined as a material whose 
diffraction patterns consist essentially of bright peaks (c.f. \cite{DefCryst}).
Since tilings made by aperiodic tiles, which we define as aperiodic tilings, are not periodic, they have served as good models for the structures of quasicrystals
when they show the diffraction patterns consisting of pure point diffraction spectrum, i.e., 
Bragg peaks only without diffuse background.
Many examples of known aperiodic tiling like Penrose tiling and Ammann tiling
show the pure point diffraction spectrum.
The objective of this article is investigating whether the
aperiodic Taylor-Socolar tiling, which is a fixed point of a substitution, 
has pure point diffraction spectrum. 

Mathematically the pure point diffraction 
spectrum is quite often studied through the
spectrum of the dynamics of the 
{\it dynamical hull}, that is,
a compact space generated by the closure of translation orbits of the tiling. 
The two notions of
pure pointedness in diffraction and dynamical spectra
are equivalent in quite a general setting \cite{LMS1, Gouere, BL,
Le-St}. 
In general tilings, almost periodicity of tilings is an equivalent criterion for the pure point spectrum. When it is restricted on substitution tilings, the almost periodicity can be easily checked by overlap coincidence. Briefly it means the following: 
when two tiles in a tiling intersect in the interior after shifting one tile by a translation of two other same type tiles in the tiling, 
one can observe a pair of same type subtiles in the same position in 
the common interior (see Subsec.\,\ref{Overlap}) 
\cite{Solomyak97, LMS2}.

An aperiodic Taylor-Socolar tiling itself does not follow tile-substitution rule strictly.
But as it is mentioned in \cite{Taylor, ST}, half-hexagonal tiles 
satisfy tile-substitution.
We should note here that the Taylor-Socolar half-hexagonal tiling is mutually locally derivable from the Taylor-Socolar tiling. 
We consider the Taylor-Socolar half-hexagonal substitution tiling
whose identical image in Taylor-Socolar tilings belongs to the dynamical hull generated by a repetitive Taylor-Socolar tiling.
We apply the substitution data on the half-hexagonal tiles of Taylor-Socolar tiles to the algorithm for checking the overlap coincidence.
The algorithm can be found in \cite{AL-web}.
As the result, we were able to check that the half-hexagonal Taylor-Socolar substitution tiling has overlap coincidence. 
So we can conclude that the aperiodic Taylor-Socolar tiling has pure point spectrum. One can also note that a dynamical hull of Taylor-Socolar tilings is invariant under the action of rotations of $\frac{n \pi}{3}$. In the diffraction pattern of a Taylor-Socolar tiling, we observed six-fold rotational symmetry.

The tiling space of Taylor-Socolar tilings with the matching rules is slightly bigger than the tiling space of  
a  Taylor-Socolar substitution tiling.
But it is shown in \cite{LM3} that the tilings in the difference have pure point spectrum, computing the total index of cosets. 
Thus any Taylor-Socolar tiling under the same matching rules have the pure point spectrum. In the case of self-similar tilings, the discrete part of the diffraction pattern, which is called Bragg spectrum, can be characterized in terms of the Fourier modules. There are three types of Bragg spectra- limit-periodic, quasiperiodic, limit-quasiperiodic. The Taylor-Socolar tilings belong to limit-periodic structure, since the expansion factor is rational. The Fourier module of this structure is an aperiodic structure which is the limit of a sequence of periodic structures (see \cite{GK}). 

 A substitution tiling with half-hexagonal shapes was known much earlier
(see \cite[Ex.10.1]{GS}). 
It is known from \cite{Fre} that the
substitution point set representing the half-hexagonal substitution tiling is a
cut-and-project set and so it has pure point spectrum. 
However Taylor-Socolar
tiling differs from \cite{Fre} in the sense
that we consider a substitution tiling reflecting the aperiodicity of
Taylor-Socolar mono-tile and have to distinguish protiles by their colors.

Various other ways to observe the pure point spectrum are pointed out in \cite{BGG, LM3}.
One can observe that there is an one-to-one almost everywhere map from a dynamical hull of Taylor-Socolar tilings to a dynamical hull of half-hexagonal substitution tilings. 
Then this induces the pure point spectrum of Taylor-Socolar tilings using the result of \cite{Fre}.
The other observation would be through checking the modular coincidence which has been introduced in \cite{LM1, LMS2}. One can see in the figures of \cite{Taylor} and \cite[Fig. 15]{ST} that $C$ or $\overline{C}$ type tiles form a sublattice structure of a whole hexagonal lattice with the expansion factor of $2$. 
It would be sufficient to check if the modular coincidence occurs with these $C$ and $\overline{C}$ type tiles. 
Furthermore \cite{LM3} provides a geometrical way to observe the limit-periodic structure in Taylor-Socolar tilings which shows that tiling can be decomposed into a superposition of periodic structure (see Fig.\,19 and Thm.\,7.1 in \cite{LM3}).

Therefore 
the pure point spectrum may not be so surprising 
in the case of Taylor-Socolar tilings. 
However comparing with the above methods, the biggest advantage 
of our method is that it is almost {\it automatic} and it can be applied to many variations of Taylor-Socolar substitutions based on hexagonal lattices, such as Penrose mono-tile tiling, with minor changes of the substitution data and investigate the diffraction spectrum of a tiling generated by it. 
Furthermore it can be applied to substitution tilings whose underlying structures are not even on lattices.


\noindent

\section{Substitution of Taylor-Socolar tiling}

\subsection{Tilings and point sets}

We briefly mention the notions of tilings and tile-substitution in $\R^2$ that we use in this paper. For more about tilings and tile-substitutions, see \cite{LP, LMS2}.

\subsubsection{Tilings}

We begin with a set of types (or colors) $\{1,\ldots,m\}$. A {\em tile} in $\R^2$ is defined as a
pair $T=(A,i)$ where $A=\supp(T)$ (the support of $T$) is a
compact set in $\R^2$, which is the closure of its interior, and
$i=l(T)\in \{1,\ldots,m\}$ is the type of $T$.
A {\em tiling} of $\R^2$ is a set $\Tk$ of tiles such that $\R^2 =
\bigcup \{\supp(T) : T \in \Tk\}$ and distinct tiles have disjoint
interiors. We always assume that any two $\Tk$-tiles with the same
color are translationally equivalent. Let $\Xi(\Tk) : = \{x \in
\R^2 : T = x + T' \ \ \mbox{for some} \ T, T' \in \Tk \}$. 
We say that a set $P$ of tiles is a {\em patch} if the number of tiles in $P$ is finite and the tiles of $P$ have mutually disjoint interiors.
We define $\Tk \cap A := \{T \in \Tk : \supp(T) \cap A \neq \emptyset\}$ for
a bounded set $A \subset \R^2$. We say that $\Tk$ is {\em
repetitive} if for every compact set $K \subset \R^2$, $\{t \in
\R^2 : \Tk \cap K = (t + \Tk) \cap K\}$ is relatively dense. 
We say that a tiling $\Tk$ has {\em finite local complexity (FLC)}
if for each radius $R > 0$ there are only finitely many
translational classes of patches whose support lies in some ball
of radius $R$.

\subsubsection{Point sets}

\noindent A {\em multi-color set} or {\em $m$-multi-color set}
in $\R^d$ is a subset $\Lb = \Lam_1 \times \dots \times \Lam_m
\subset \R^d \times \dots \times \R^d$ \; ($m$ copies) where
$\Lam_i \subset \R^d$. We also write $\Lb = (\Lam_1, \dots,
\Lam_m) = (\Lam_i)_{i\le m}$. Recall that a Delone set is a
relatively dense and uniformly discrete subset of $\R^d$. We say
that $\Lb=(\Lambda_i)_{i\le m}$ is a {\em Delone multi-color set}
in $\R^d$ if each $\Lambda_i$ is Delone and
$\supp(\Lb):=\cup_{i=1}^m \Lambda_i \subset \R^d$ is Delone. 
We say that $\Lam \subset \R^d$ is a {\em Meyer set}
if it is a Delone set and $\Lam - \Lam$ is uniformly discrete (\cite{Lag}).
The types (or
colors) of points on Delone multi-color sets have the same
concept as the colors of tiles on tilings.

\subsection{Tile substitution and associated substitution Delone multi-color set}

We say that a linear map $Q : \R^2 \rightarrow \R^2$ is {\em
expansive} if all the eigenvalues of $Q$
 lie outside the closed unit disk in $\C$.

\begin{defi}\label{def-subst}
{\em Let $\Ak = \{T_1,\ldots,T_m\}$ be a finite set of tiles in
$\R^2$ such that $T_i=(A_i,i)$; we will call them {\em
prototiles}. Denote by ${\mathcal{P}}_{\Ak}$ the set of non empty
patches. 
We say that $\Omega: \Ak \to {\mathcal{P}}_{\Ak}$ is a {\em
tile-substitution} (or simply {\em substitution}) with an
expansive map $Q$ if there exist finite sets $\Dk_{ij}\subset
\R^2$ for $i,j \le m$ such that
\begin{equation}
\Omega(T_j)=
\{u+T_i:\ u\in \Dk_{ij},\ i=1,\ldots,m\}
\label{subdiv}
\end{equation}
with
\begin{eqnarray} \label{tile-subdiv}
Q A_j = \bigcup_{i=1}^m (\Dk_{ij}+A_i) \ \ \  \mbox{for} \  j\le m.
\end{eqnarray}
Here all sets in the right-hand side must have disjoint interiors;
it is possible for some of the $\Dk_{ij}$ to be empty.}
\end{defi}

\noindent The substitution (\ref{subdiv}) is extended to all
translates of prototiles by $\Omega(x+T_j)= Q x + \Omega(T_j)$ and
to patches and tilings by $\Omega(P)=\bigcup\{\Omega(T):\ T\in
P\}$. The substitution $\Omega$ can be iterated, producing larger
and larger patches $\Omega^k(P)$.
We say that $\Tk$ is a {\em substitution tiling} if $\Tk$ is a
tiling and $\Omega(\Tk) = \Tk$ with some substitution $\Omega$. In
this case, we also say that $\Tk$ is a {\em fixed point} of
$\Omega$. We say that a substitution tiling is {\em primitive} if
the corresponding substitution matrix $S$, with $S_{ij}= \sharp
(\Dk_{ij})$, is primitive. A repetitive fixed point of a primitive
tile-substitution with FLC is called a {\em self-affine tiling}.
If $\Tk = \lim_{n \to \infty} \Omega^n(P)$, we say that $P$ is a {\em generating patch}.

We say that $\Lb=(\Lambda_i)_{i\le m}$ is a {\em Delone
multi-color set} in $\R^2$ if each $\Lambda_i$ is Delone and
$\supp(\Lb):=\cup_{i=1}^m \Lambda_i \subset \R^2$ is Delone.
Any tiling $\Tk$ can be converted into a Delone multi-color set by simply
choosing a point $x_{(A,i)}$ for each tile $(A,i)$ so that the
chosen points for tiles of the same type are in the same relative
position in the tile: $x_{(g+A,i)}= g + x_{(A,i)}$. We define
$\Lambda_i := \{x_{(A,i)}:\ (A,i) \in \Tk \}$ and $\Lb :=
(\Lambda_i)_{i \le m}$. Clearly $\Tk$ can be reconstructed from
$\Lb$ given the information about how the points lie in their
respective tiles. This bijection establishes a topological
conjugacy of $(X_{\Lbs},\R^2)$ and $(X_{\Tk},\R^2)$. Concepts and
theorems can clearly be interpreted in either language (FLC, repetitivity, pure point dynamical spectrum, etc.).

If a self-affine tiling $\Tk = \{T_j+ \Lam_j : j \le m \}$ is
given, we get an associated substitution Delone multi-color set
$\Lb_{\Tk} = (\Lam_i)_{i \le m}$ of $\Tk$ (see
\cite[Lemma\,5.4]{Lee}).

\subsection{Two equivalent criteria for pure point spectrum} \label{Two-pure-pointedness}
There are two notions of pure pointness in the study of tilings - pure point dynamical spectrum and pure point diffraction spectrum. We briefly give the definitions of them.

Let $\Tk$ be a tiling in $\R^2$. We define the space of tilings as
the orbit closure of $\Tk$ under the translation action: $X_{\Tk}
= \overline{\{-h + \Tk : h \in \R^2 \}}$, in the well-known
``local topology'': for a small $\epsilon
> 0$ two point sets $\Sk_1, \Sk_2$ are $\epsilon$-close if $\Sk_1$
and $\Sk_2$ agree on the ball of radius $\epsilon^{-1}$ around the
origin, after a translation of size less than $\epsilon$. The
group $\R^2$ acts on $X_{\Tk}$ by translations which are obviously
homeomorphisms, and we get a topological dynamical system
$(X_{\Tk},\R^2)$. Let $\mu$ be an ergodic invariant Borel
probability measure for the dynamical system $(X_{\Tk},\R^2)$. We
consider the associated group of unitary operators $\{U_g\}_{g\in
\R^2}$ on $L^2(X_{\Tk},\mu)$ for which $U_g f(\Sk) = f(-g+\Sk)$.
The dynamical system $(X_{\Tk},\mu,\R^2)$ is said
to have {\em pure point(or pure discrete) spectrum} if the linear
span of the eigenfunctions is dense in $L^2(X_{\Tk}, \mu)$.

On the other hand, there is a notion of pure point diffraction
spectrum which characterizes quasicrystals. Let $\Lb = (\Lam_i)_{i
\le m}$ be a multi-color point set in $\R^2$. We consider a
measure of the form $\nu=\sum_{i \le m} a_i \delta_{\Lam_i}$,
where $\delta_{\Lam_i} = \sum_{x \in \Lam_i} \delta_x$ and $a_i
\in \C$. The autocorrelation of $\nu$ is
\[ \gamma(\nu) = \lim_{n \to \infty} \frac{1}{\Vol(B_n)} (\nu|_{B_n} \ast \widetilde{\nu}|_{B_n}),\]
where $\nu|_{B_n}$ is a measure of $\nu$ restricted on the ball
$B_n$ of radius $n$ and $\widetilde{\nu}$ is the measure, defined
by $\widetilde{\nu}(f) = \overline{\nu(\widetilde{f})}$, where $f$
is a continuous function with compact support and
$\widetilde{f}(x) = \overline{f(-x)}$. The diffraction measure of
$\nu$ is the Fourier transform $\widehat{\gamma(\nu)}$ of the
autocorrelation (see \cite{Hof}). When the diffraction measure
$\widehat{\gamma(\nu)}$ is a pure point measure, we say that $\Lb$
has pure point diffraction spectrum and so $\Lb$ has the
structures of quasicrystals.

It turns out that the two notions of pure pointedness are same, i.e,
the pure point dynamical spectrum of $(X_{\Tk}, \R^2, \mu)$ is
equivalent to the pure point diffractivity of $\Lb_{\Tk}$
\cite{LMS1, Gouere, BL, Le-St}.

\subsection{Overlap} \label{Overlap}

A triple $(u, y, v)$, with $u + T_i, v + T_j \in \Tk$ and $y \in
\Xi(\Tk)$, is called an {\em overlap} if
\[ (u+A_i-y)^{\circ} \cap (v+A_j)^{\circ} \neq \emptyset, \]
where $A_i = \supp(T_i)$ and $A_j = \supp(T_j)$. We define
$(u+A_i-y) \cap (v+A_j)$ the {\em support of an overlap} $(u, y,
v)$ and denote it by $\supp(u,y,v)$.
An overlap $(u, y, v)$ is a {\em coincidence} if
\[ \mbox{$u-y = v$ and $u + T_i, v + T_i \in \Tk$ for some $i \le m$}.\]
Let $\mathcal{O} = (u, y, v)$ be an overlap in $\Tk$, we define
{\em $k$-th inflated overlap}
\begin{eqnarray*}
 {\Phi}^k \mathcal{O} = \{(u', Q^k y, v') \, :
  u' \in \Phi^k(u), v' \in \Phi^k(v), \ \mbox{and $(u',Q^ky,v')$ is an overlap} \}.
\end{eqnarray*}

\begin{defi} \label{def-overlapCoincidence}
{\em We say that a self-affine tiling $\Tk$ admits an {\em overlap
coincidence} if there exists $\ell \in \Z_+$ such that for each
overlap $\mathcal{O}$ in $\Tk$, ${\Phi}^{\ell} \mathcal{O}$
contains a coincidence.}
\end{defi}

When $\Tk$ is a self-affine tiling in $\R^2$ such that $\Xi(\Tk)$
is a Meyer set, $(X_{\Tk}, \R^2, \mu)$ has a pure point
dynamical spectrum if and only if $\Tk$ admits an overlap
coincidence\cite{LMS2, Lee}. So we will check the pure point spectrum of Taylor-Socolar tilings through the computation of overlap coincidence in the next subsection.

\subsection{Taylor-Socolar half-hexagonal substitution tiling}

In the half-hexagonal substitution  tiling, there are 14 half-hexagonal
prototiles which come from dividing 7 hexagonal prototiles $A, B, C, D, E, F, G$ into
the left and the right (see \cite[Fig. 15]{ST}). Since the substitution tiling we defined in (\ref{def-subst}) requires finite prototiles up to only translations, we
need to treat the rotated types and reflected types of the 14 half-hexagonal
prototiles as different prototiles. So we consider a substitution
tiling with 168 prototiles. Using the algorithm in
\cite{AL-web} which is originated from \cite{AL} and made for the computation of the Taylor-Socolar tilings, we check if the half-hexagonal substitution Taylor-Socolar tiling has pure
point spectrum.

Now a question would be ``when the dynamics of the half-hexagonal
substitution tiling has pure point spectrum, can we infer that the dynamics of
the original hexagonal tiling also has pure point spectrum?''.
Let $\Tk$ be a fixed point of a primitive substitution and $\Lb_{\Tk} = (\Lambda_i)_{i \le m}$. It is shown in \cite[Lemma A.6.]{LMS2} that $\Tk$ has uniform cluster frequencies(UCF), i.e. for any $\Tk$-patch $P$, there exists 
\[ \freq(P, \Tk) := \lim_{n \to \infty} \frac{L_P(h+ B_n)}{\mbox{Vol}(B_n)}\]
uniformly in $h \in \R^2$.
From \cite[Theorem 3.2]{LMS1}, the measure $\nu = \sum_{i \le m} a_i
    \delta_{\Lambda_i}$ has pure point diffraction spectrum,
    for any choice of complex numbers $(a_i)_{i\le m}$ if and only if 
each measure $\delta_{\Lambda_i}$ has pure
    point diffraction spectrum.
By the construction of the half-hexagonal Taylor-Socolar tiling, we can take a
substitution point set $\Lb$ representing the half-hexagonal Taylor-Socolar tiling 
to include a substitution point set $\Gb$ representing the
original hexagonal Taylor-Socolar tiling. So from \cite[Theorem 3.2]{LMS1},
we can conclude that Taylor-Socolar tiling has pure point diffraction spectrum.

\begin{figure}
\begin{center}
\includegraphics[width=0.3\textwidth]{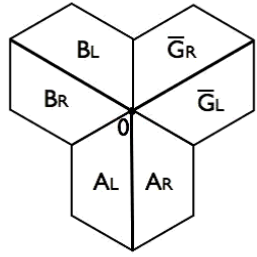}
\caption{A generating patch around the origin of the substitution tiling with half-hexagonal Taylor-Socolar tiles} \label{TS-generating-patch}
\end{center}
\end{figure}

As a generating patch, one can start with the patch as shown in Figure \ref{TS-generating-patch}.
Since this patch is contained in the next inflated patch after the substitution, it gives a fixed tiling under the substitution.

\begin{figure}
\begin{center}
\includegraphics[width=0.9\textwidth]{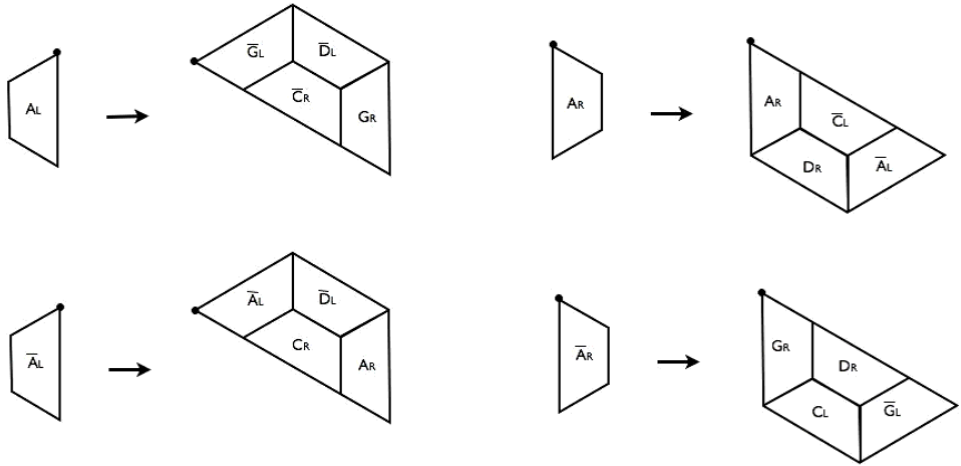}
\caption{Substitution of half-hexagonal Taylor-Socolar tile types $A_L$, $\overline{A}_L$, $A_R$, and $\overline{A}_R$. The black dot indicates the origin.} \label{TS-substitution}
\end{center}
\end{figure}

For the computational reason, we give a tile-substitution whose expansion involves rotation and reflection.
It is possible to use the second iteration of half-hexagons without rotation and reflection parts on expansions. But in this case, the substitution gets bigger.  
The tile-substitution for $A_L, A_R, \overline{A}_L, \overline{A}_R$ is shown in Figure \ref{TS-substitution}. For other half-hexagonal tiles, the figures of tile-substitution is similar. We give the precise tile-substitution below. One can check the computational algorithm in \cite{AL-web}.
Let
\[Q = 2 \cdot R_{ot} \cdot R_{ef} =
2 \left( \begin{array}{rr}
         \cos \frac{\pi}{3} & -\sin \frac{\pi}{3} \\
         \sin \frac{\pi}{3} & \cos \frac{\pi}{3}
         \end{array}
 \right) \left( \begin{array}{rr}
         -1 & 0 \\
          0 & 1
         \end{array}
 \right)
= 2 \left( \begin{array}{rr}
         -\cos \frac{\pi}{3} & -\sin \frac{\pi}{3} \\
         -\sin \frac{\pi}{3} & \cos \frac{\pi}{3}
         \end{array}
 \right),\]
where $R_{ot}$ is a rotation of $\frac{\pi}{3}$ counter clockwise through the origin
and $R_{ef}$ is a reflection through $y$-axis. The tile
substitution is given as follows; let $\omega = \left(
        \begin{array}{rr}
         \cos \frac{\pi}{3} & -\sin \frac{\pi}{3} \\
         \sin \frac{\pi}{3} & \cos \frac{\pi}{3}
         \end{array}
 \right)$ and $u = \left(
        \begin{array}{r}
         \cos \frac{\pi}{6}  \\
         \sin \frac{\pi}{6}
         \end{array}
 \right)$. We denote
$(S_X)_n := \omega^n S_X $, where $m \in \Z$, $S \in \{A, B, C, D,
E, F, G, \overline{A}, \overline{B}, \overline{C}, \overline{D},
\overline{E}, \overline{F}, \overline{G} \}$, and $X \in \{L,
R\}$. 

The tile-substitution rule for $A_L$ and $A_R$ is the following.
For $0 \le n \le 5$,
\begin{eqnarray*}
Q (A_L)_n & = & (\overline{G}_L)_{2-n} \cup ((\overline{D}_L)_{1-n} + 2 u) \cup (
(\overline{C}_R)_{1-n} + \omega^5 u) \cup ((G_R)_{3-n} + 4 \omega^5 u) \\
 Q (A_R)_n & = & (A_R)_{-n} \cup ((D_R)_{1-n} + 2 \omega ^4 u) \cup ((\overline{C}_L)_{1-n} +
\omega^5 u) \cup ((\overline{A}_L)_{5-n} + 4 \omega^5 u) \\
Q (\overline{A}_L)_n & = & (\overline{A}_L)_{2-n} \cup ((\overline{D}_L)_{1-n} + 2 u) \cup
((C_R)_{1-n} + \omega^5 u) \cup ((A_R)_{3-n} + 4 \omega^5 u) \\
 Q (\overline{A}_R)_n & = & (G_R)_{-n} \cup ((D_R)_{1-n} + 2 \omega ^4 u) \cup ((C_L)_{1-n} +
\omega^5 u) \cup ((\overline{G}_L)_{5-n} + 4 \omega^5 u). 
\end{eqnarray*}

The tile substitution rule for other types of tiles are similar. We give the rules for the convenience to the readers.
{\footnotesize
\begin{eqnarray*}
Q (B_L)_n & = & (B_L)_{2-n} \cup ((\overline{F}_L)_{1-n} + 2 u)
\cup
((\overline{C}_R)_{1-n} + \omega^5 u) \cup ((G_R)_{3-n} + 4 \omega^5 u) \\
 Q (B_R)_n & = & (\overline{G}_R)_{-n} \cup ((F_R)_{1-n} + 2 \omega^4 u) \cup
 ((\overline{C}_L)_{1-n} +
\omega^5 u) \cup ((\overline{A}_L)_{5-n} + 4 \omega^5 u) \\
Q (\overline{B}_L)_n & = & (G_L)_{2-n} \cup ((\overline{F}_L)_{1-n} + 2 u) \cup
((C_R)_{1-n} + \omega^5 u) \cup ((A_R)_{3-n} + 4 \omega^5 u) \\
 Q (\overline{B}_R)_n & = & (\overline{B}_R)_{-n} \cup ((F_R)_{1-n} + 2 \omega^4 u) \cup ((C_L)_{1-n} +
\omega^5 u) \cup ((\overline{G}_L)_{5-n} + 4 \omega^5 u),
\end{eqnarray*}
\begin{eqnarray*}
Q (C_L)_n & = & (F_L)_{2-n} \cup ((E_L)_{1-n} + 2 u) \cup
((\overline{C}_R)_{1-n} + \omega^5 u) \cup ((\overline{F}_R)_{3-n} + 4 \omega^5 u) \\
 Q (C_R)_n & = & (\overline{D}_R)_{-n} \cup ((\overline{E}_R)_{1-n} + 2 \omega^4 u) \cup ((\overline{C}_L)_{1-n} +
\omega^5 u) \cup ((D_L)_{5-n} + 4 \omega^5 u) \\
Q (\overline{C}_L)_n & = & (D_L)_{2-n} \cup ((E_L)_{1-n} + 2 u) \cup
((C_R)_{1-n} + \omega^5 u) \cup ((\overline{D}_R)_{3-n} + 4 \omega^5 u) \\
 Q (\overline{C}_R)_n & = & (\overline{F}_R)_{-n} \cup ((\overline{E}_R)_{1-n} + 2 \omega^4 u) \cup
 ((C_L)_{1-n} + \omega^5 u) \cup ((F_L)_{5-n} + 4 \omega^5 u),
\end{eqnarray*}
\begin{eqnarray*}
Q (D_L)_n & = & (\overline{B}_L)_{2-n} \cup ((\overline{D}_L)_{1-n} + 2 u) \cup
((\overline{C}_R)_{1-n} + \omega^5 u) \cup ((B_R)_{3-n} + 4 \omega^5 u) \\
 Q (D_R)_n & = & (\overline{A}_R)_{-n} \cup ((E_R)_{1-n} + 2 \omega^4 u) \cup
 ((\overline{C}_L)_{1-n} +
\omega^5 u) \cup ((A_L)_{5-n} + 4 \omega^5 u) \\
Q (\overline{D}_L)_n & = & (A_L)_{2-n} \cup ((\overline{E}_L)_{1-n} + 2 u) \cup
((C_R)_{1-n} + \omega^5 u) \cup ((\overline{A}_R)_{3-n} + 4 \omega^5 u) \\
 Q (\overline{D}_R)_n & = & (B_R)_{-n} \cup ((D_R)_{1-n} + 2 \omega^4 u) \cup
 ((C_L)_{1-n} + \omega^5 u) \cup ((\overline{B}_L)_{5-n} + 4 \omega^5 u),
\end{eqnarray*}
\begin{eqnarray*}
Q (E_L)_n & = & (B_L)_{2-n} \cup ((\overline{E}_L)_{1-n} + 2 u) \cup
((\overline{C}_R)_{1-n} + \omega^5 u) \cup ((\overline{B}_R)_{3-n} + 4 \omega^5 u) \\
 Q (E_R)_n & = & (\overline{G}_R)_{-n} \cup ((E_R)_{1-n} + 2 \omega^4 u) \cup
 ((\overline{C}_L)_{1-n} +
\omega^5 u) \cup ((G_L)_{5-n} + 4 \omega^5 u) \\
Q (\overline{E}_L)_n & = & (G_L)_{2-n} \cup ((\overline{E}_L)_{1-n} + 2 u) \cup
((C_R)_{1-n} + \omega^5 u) \cup ((\overline{G}_R)_{3-n} + 4 \omega^5 u) \\
 Q (\overline{E}_R)_n & = & (\overline{B}_R)_{-n} \cup ((E_R)_{1-n} + 2 \omega^4 u) \cup ((C_L)_{1-n} +
\omega^5 u) \cup ((B_L)_{5-n} + 4 \omega^5 u),
\end{eqnarray*}
\begin{eqnarray*}
Q (F_L)_n & = & (B_L)_{2-n} \cup ((\overline{F}_L)_{1-n} + 2 u) \cup
((\overline{C}_R)_{1-n} + \omega^5 u) \cup ((B_R)_{3-n} + 4 \omega^5 u) \\
 Q (F_R)_n & = & (\overline{G}_R)_{-n} \cup ((E_R)_{1-n} + 2 \omega^4 u) \cup
 ((\overline{C}_L)_{1-n} +
\omega^5 u) \cup ((A_L)_{5-n} + 4 \omega^5 u) \\
Q (\overline{F}_L)_n & = & (G_L)_{2-n} \cup ((\overline{E}_L)_{1-n} + 2 u) \cup
((C_R)_{1-n} + \omega^5 u) \cup ((\overline{A}_R)_{3-n} + 4 \omega^5 u) \\
 Q (\overline{F}_R)_n & = & (\overline{B}_R)_{-n} \cup ((F_R)_{1-n} + 2 \omega^4 u) \cup ((C_L)_{1-n} +
\omega^5 u) \cup ((\overline{B}_L)_{5-n} + 4 \omega^5 u),
\end{eqnarray*}
\begin{eqnarray*}
Q (G_L)_n & = & (\overline{B}_L)_{2-n} \cup ((\overline{D}_L)_{1-n} + 2 u) \cup
((\overline{C}_R)_{1-n} + \omega^5 u) \cup ((G_R)_{3-n} + 4 \omega^5 u) \\
 Q (G_R)_n & = & (\overline{A}_R)_{-n} \cup ((F_R)_{1-n} + 2 \omega^4 u) \cup
 ((\overline{C}_L)_{1-n} +
\omega^5 u) \cup ((\overline{A}_L)_{5-n} + 4 \omega^5 u) \\
Q (\overline{G}_L)_n & = & (A_L)_{2-n} \cup ((\overline{F}_L)_{1-n} + 2 u) \cup
((C_R)_{1-n} + \omega^5 u) \cup ((A_R)_{3-n} + 4 \omega^5 u) \\
 Q (\overline{G}_R)_n & = & (B_R)_{-n} \cup ((D_R)_{1-n} + 2 \omega^4 u) \cup
 ((C_L)_{1-n} + \omega^5 u) \cup ((\overline{G}_L)_{5-n} + 4 \omega^5 u).
\end{eqnarray*}
}


Let us describe how the algorithm works in general term. When tile-substitution date $\Omega$ is given, we consider a substitution Delone multi-color set $\Lb$ which is fixed under the substitution. To build the set $\Lb$, we need to find a point $x \in \Lb$ which is fixed under the substitution. Applying the substitutions to $\{x\}$ infinitly many times, we can easily obtain a point set $\Lb$ which is fixed under the substitution. It is sufficient to check the overlap coincidence for all the overlaps which occur by finite translation vectors of same type tiles in the tiling. From the Meyer property, the number of overlaps are finite. After collecting all the overlaps, we can check overlap coincidence for each overlap applying the substitution many times. Here the number of times of applying the substitution can be limited by the number of overlaps. So the algorithm will be terminated. The detail is given in \cite{AL}.

The computation of overlap coincidence of Taylor-Socolar tiling takes rather long time comparing to other examples in \cite{AL}. We guess that it is due to the number of prototiles (168) which is much more than other examples. 
It is an interesting question whether the complexity of the substitution rule is 
related with the computation time of the algorithm when the number of prototiles is same.

\section{Acknowledgment}

\noindent
We are grateful to Michael Baake and Franz G\"ahler for their interest and helpful
suggestions.

\vspace{6mm}

{\footnotesize
a: 
Institute of Mathematics, University of Tsukuba, 1-1-1 Tennodai, \\
 \hspace*{2.5em} Tsukuba, Ibaraki, Japan (zip:350-0006);
\email{akiyama@math.tsukuba.ac.jp}

\smallskip
b: Dept. of Math. Edu., Kwandong University, 24, 579 Beon-gil, Beomil-ro, Gangneung,  \\  
\hspace*{2.5em} Gangwon-do, 210-701 Republic of Korea;
\email{jylee@kd.ac.kr,  jeongyuplee@yahoo.co.kr} }

\end{document}